\documentclass{commat}

\usepackage{caption}
\usepackage{lscape}
\usepackage{longtable}
\usepackage{pdflscape}
\usepackage{graphicx}
\usepackage{enumerate}
\usepackage{xcolor}
\usepackage[utf8]{inputenc}

\DeclareMathOperator{\Bro}{Bro}
\newenvironment{items}{\begin{list}{$\alph{item})$}{\labelwidth18pt \leftmargin28pt \topsep3pt \itemsep1pt \parsep0pt}}{\end{list}}

\title{%
   Computing subalgebras and $\mathbb{Z}_2$-gradings 
        of simple Lie algebras over finite fields
    }

\author{%
    Bettina Eick and Tobias Moede
    }

\affiliation{
    \address{
    Institute for Analysis and Algebra, TU Braunschweig, Germany
    }
    \email{%
    beick@tu-bs.de, t.moede@tu-bs.de
    }
}

\abstract{%
    This paper introduces two new algorithms for Lie algebras over finite fields and applies them to the investigate the known simple Lie algebras of dimension at most $20$ over the field $\mathbb{F}_2$ with two elements. The first algorithm is a new approach towards the construction of $\mathbb{Z}_2$-gradings of a Lie algebra over a finite field of characteristic $2$. Using this, we observe that each of the known simple Lie algebras of dimension at most $20$ over $\mathbb{F}_2$ has a $\mathbb{Z}_2$-grading and we determine the associated simple Lie superalgebras. The second algorithm allows us to compute all subalgebras of a Lie algebra over a finite field. We apply this to compute the subalgebras, the maximal subalgebras and the simple subquotients of the known simple Lie algebras of dimension at most $16$ over $\mathbb{F}_2$ (with the exception of the $15$-dimensional Zassenhaus algebra).
    }



\VOLUME{30}
\NUMBER{2}
\firstpage{37}
\DOI{https://doi.org/10.46298/cm.10193}

\begin{paper}

\section{Introduction}
\label{sec:introduction}

The classification of the finite-dimensional simple Lie algebras depends
heavily on their underlying field. For algebraically closed fields of
characteristic zero the classification has been achieved long ago and
is folklore nowadays. For algebraically closed fields of characteristic
$p \geq 5$ the classification has been completed more recently, see
\cite{PSt06}. A classification over fields of characteristic $p \in 
\{2,3\}$ has not been achieved so far. The available evidence suggests 
that the classification over fields of characteristic $2$ will differ significantly 
from the other cases. This motivates the computational investigation of the 
known simple Lie algebras over the field with two elements with the aim to 
gain further insight into their structure.

In this paper we introduce two algorithms to investigate Lie algebras
over finite fields. Our first method determines $\mathbb{Z}_2$-gradings.
A well-known construction for $\mathbb{Z}_2$-gradings of Lie algebras over fields 
$\mathbb{F}$ of characteristic different from $2$ uses the Cartan decomposition of $L$: If
$\theta \colon L \rightarrow L$ is an automorphism with $\theta^2 = 1$, then $\theta$
is diagonalizable with eigenvalues $1$ and $-1$ and its eigenspaces 
form a $\mathbb{Z}_2$-grading. We introduce an alternative construction for fields
of characteristic $2$: If $l \in L$ with $(ad_L(l))^2 = ad_L(l)$, 
then $ad_L(l)$ is diagonalizable with eigenvalues $0$ and $1$ and its
eigenspaces form a $\mathbb{Z}_2$-grading. We call $l$ an {\em
idempotent} and we introduce an algorithm to compute the idempotents in a 
Lie algebra $L$ over a finite field $\mathbb{F}$ with
$char(\mathbb{F}) = 2$. We use this algorithm to determine the $Aut(L)$-orbits
of idempotents in the known simple Lie algebras of dimension at most
$20$ over the field $\mathbb{F}_2$ with two elements, see Section \ref{results}
for details. Based on this computation, we propose the following conjecture.

\begin{conjecture}
Every finite-dimensional simple Lie algebra over a field of characteristic $2$ has a non-central idempotent, and hence a non-degenerate $\mathbb{Z}_2$-grading.
\end{conjecture}

Lie algebras with $\mathbb{Z}_2$-grading are closely related to Lie superalgebras.
The latter play a role in physics where they are used to describe the 
mathematics of supersymmetry, we refer to \cite{Kac77}, \cite{BeBou18} 
and \cite{KrLe18} for details on Lie superalgebras. In \cite{BoLe14} it
has been shown that each $\mathbb{Z}_2$-grading of a simple Lie algebra over a
field of characteristic $2$ determines a simple Lie superalgebra. We 
recall this construction in Theorem \ref{super} below for completeness, 
and use it to determine the simple Lie 
superalgebras arising from the
known simple Lie algebras of dimension at most $20$ over $\mathbb{F}_2$. See
Section \ref{results} for details.

Our second aim in this paper is to describe an algorithm that computes all 
subalgebras of a finite Lie algebra $L$ up to the action of $Aut(L)$. 
We apply this algorithm to determine all subalgebras of the currently 
known simple Lie algebras of dimension at most $16$ over the field $\mathbb{F}_2$ 
with two elements, except for the $15$-dimensional Zassenhaus algebra, which 
we denote by $W(4)$. The Lie algebra $W(4)$ has
more than two million orbits of subalgebras under the action of its
automorphism group and it is the only case in our considered range 
that we could not complete. Our algorithm also allows to determine 
maximal subalgebras, and the simple subquotients of a Lie algebra over
a finite field. We exhibit our computational results in Section 
\ref{results} below.

Our algorithms are implemented in the computer algebra system GAP 
\cite{GAP}. We use the FinLie package \cite{finlie} for the computation
of automorphism groups of Lie algebras over $\mathbb{F}_2$ and the list 
of known low-dimensional simple Lie algebras over the field $\mathbb{F}_2$, which is
also a part of FinLie. This list contains the complete classification of Lie
algebras of dimension at most $9$ obtained by Vaughan-Lee \cite{VL06}, 
the list determined by Eick \cite{Eic10} using computational methods,
three Lie algebras that have been computed recently using the methods 
of \cite{Eic10}, and a Lie algebra described by Skryabin \cite{Skr98}. 
Generators for the three new Lie algebras can be found in Appendix 
\ref{appendix} below. 
\medskip

\section{Idempotents, gradings and Lie superalgebras}

Let $L$ be a Lie algebra and let $B$ be an abelian additive group. The Lie algebra 
$L$ is $B$-graded if $L$ can be written as a direct sum of vector spaces
\[ L = \bigoplus_{b \in B} L_b \]
so that $[L_a, L_b] \subseteq L_{a+b}$ holds for all $a,b \in B$. Such a
grading is also called a {\em group grading}. We say that a group grading 
is {\em degenerate} if there exists $b \in B$ with $L_b = \{0\}$.  Patera \& 
Zassenhaus \cite{PZa89} initiated a first systematic study of arbitrary 
non-degenerate gradings. We also refer to \cite{Cal14} and \cite{ElKo13} 
for further background on gradings.

Of particular interest are gradings with $B \cong \mathbb{Z}_2$ being cyclic of order 
$2$. In the following section we introduce a construction for such
a type of grading.

\subsection{Idempotents and gradings}

Let $L$ be a Lie algebra over an arbitrary field $\mathbb{F}$. For $l \in L$ and
$a \in \mathbb{F}$ we define the vector space $E_a(l) = \{ h \in L \mid [l,h] = ah\}$ to be the eigenspace
of $ad_L(l)$ to the eigenvalue $a$ and 
\[ E(l) = \bigoplus \limits_{a \in \mathbb{F}} E_a(l) .\]
An element $l \in 
L$ is {\em diagonalizable} if $ad_L(l)$ is diagonalizable over $\mathbb{F}$. 
We denote by $V(l) = \langle a \in \mathbb{F} \mid E_a(l) \neq \{0\} \rangle$ the additive
subgroup of $\mathbb{F}$ generated by the eigenvalues of $ad_L(l)$. 

\begin{lemma}
\label{mult}
Let $L$ be a Lie algebra over an arbitrary field and let $l \in L$. 
\begin{items}
\item[\rm (a)]
$E(l)$ is the direct sum of its non-zero subspaces $E_a(l)$ and 
$[E_a(l), E_b(l)] \subseteq E_{a+b}(l)$ holds for all $a, b \in \mathbb{F}$.
\item[\rm (b)]
Suppose that $l$ is diagonalizable. Then $E(l)$ is a group grading 
of $L$ for the group $V(l)$.
\end{items}
\end{lemma}

\begin{proof}
(a)
The set $E_a(l)$ is the eigenspace to the eigenvalue $a$ of $ad_L(l)$.
Hence $E_a(l)$ is a subspace and the eigenspaces form a direct sum.
Now let $h \in E_a(l)$ and $k \in E_b(l)$. Then $[l, [ h,k]] 
= [[l,h],k]+[h,[l,k]] = a [h,k] + b [h,k] 
= (a+b) [h,k]$. Hence $[h,k] \in E_{a+b}(l)$.

(b) Follows readily from (a).
\end{proof}

Lemma \ref{mult} leads to the following central observation.

\begin{theorem} \label{idempotentgrading}
Let $L$ be a field of characteristic $2$, and let $l$ be a non-central
idempotent of $L$. Then $L=L_0 \oplus L_1$ is a non-degenerate $\mathbb{Z}_2$-grading, 
where $L_0 = E_0(l)$ and $L_1 = E_1(l)$.
\end{theorem}

\begin{proof}
Let $m = ad_L(l)$. As $m^2 = m$, it follows that $m(m-1)=0$ and thus 
the minimal polynomial of $m$ divides $x(x-1)$. This implies that $m$
is diagonalizable with $L_0 \oplus L_1 = L$. Next $L_0 \neq \{0\}$,
since $l \in L_0$ and $L_1 \neq \{0\}$, since $l \not \in Z(L)$.
In summary, we obtain a non-degenerate $\mathbb{Z}_2$-grading for $L$.
\end{proof}

\begin{remark}
Let $\alpha \in Aut(L)$ and let $l_1, l_2 \in L \setminus Z(L)$ be idempotents
with $l_2 = \alpha(l_1)$. If $L=L_0 \oplus L_1$ is the grading associated
to $l_1$ as in Theorem \ref{idempotentgrading}, then $\alpha(L_0) \oplus \alpha(L_1)$
is the grading associated to $l_2$, i.e., both gradings belong to the same
$Aut(L)$-orbit.
\end{remark}

The idempotents of a Lie algebra $L$ over a field $\mathbb{F}$ with $char(\mathbb{F})=2$
can be computed easily. Let $b_1,\ldots,b_n$ be a basis of $L$ and let 
$x_1, \ldots, x_n$ be indeterminates. If $x = x_1 b_1 + \ldots + x_n b_n$, 
then $ad_L(x) = x_1 ad_L(b_1) + \ldots + x_n ad_L(b_n)$. Hence
\[ ad_L(x)^2 - ad_L(x) = 
   \sum_{i=1}^n \sum_{j=1}^n x_i x_j ad_L(b_i) ad_L(b_j) 
   - \sum_{i=1}^n x_i ad_L(b_i).\]
The equation $ad_L(x)^2-ad_L(x) = 0$ now translates to $n^2$ polynomial
equations in the indeterminates $x_1, \ldots, x_n$ over the field $\mathbb{F}$.
We obtain the idempotents of $L$ by determining all $(x_1, \ldots, x_n)
\in \mathbb{F}^n$ solving all of the $n^2$ equations. 

\subsection{Simple Lie superalgebras from simple Lie algebras}

We first recall the definition of a Lie superalgebra over a field $\mathbb{F}$ of
characteristic $2$ from \cite{KrLe18}, see also \cite{BoLe14} and \cite{BeBou18}. 
Let $L$ be a Lie algebra over a field $\mathbb{F}$ of characteristic $2$. Then $L$ 
is a Lie superalgebra if there exists a $\mathbb{Z}_2$-grading $L = L_0 \oplus L_1$ 
and a map $s$ (called {\em squaring})

\[ s \colon L_1 \rightarrow L_0,\; x \mapsto s(x), \]

such that $s(\alpha x) = \alpha^2 s(x)$ for all $\alpha \in \mathbb{F}, x \in L_1$, 
and  $[x,y] = s(x+y) - s(x) - s(y)$ for all $x,y \in L_1$, and $[s(x), y] = [x, [x,y]]$ for all 
$x \in L_1, y \in L$. The later condition translates to $ad_L(s(x)) = 
(ad_L(x))^2$ for $x \in L_1$, and, if the center of $L$ is trivial, then this
allows to determine $s(x)$ for $x \in L_1$ if it exists. The following 
is in part also proved in \cite{BoLe14}. We include a proof here for 
completeness.

\begin{theorem}
\label{super}
Let $L$ be a finite-dimensional $\mathbb{Z}_2$-graded Lie algebra with trivial center
over a field $\mathbb{F}$ of characteristic $2$. Then $L$ embeds into a 
finite-dimensional Lie superalgebra $S$ over $\mathbb{F}$ such that
\begin{items}
\item[\rm (a)]
as a Lie algebra, $L$ is an ideal of $S$ and $S/L$ is abelian of dimension
at most $\dim(L_1)$. 
\item[\rm (b)]
if $L$ is a simple Lie algebra, then $S$ is a simple Lie superalgebra.
\end{items}
\end{theorem}

\begin{proof}
Let $\varphi \colon L \rightarrow End(L),\; l \mapsto ad_L(l)$. As the center of $L$ is trivial,
the homomorphism $\varphi$ is injective and $L \cong \varphi(L)$. Therefore, in the 
following we identify $L$ with $\varphi(L)$.
\medskip

As a first step, we recall some details of the arithmetic of $End(L)$.
The Lie bracket in $End(L)$ is given by $[x,y] = xy -yx = xy+ yx$, since 
$char(\mathbb{F}) = 2$. For $x,y \in End(L)$ and $a,b \in \mathbb{F}$ we obtain that 

\begin{align}
(ax+by)^2 &= (ax)^2 + axby + byax + (by)^2 \label{eq1} \\
          &= a^2x^2 + b^2y^2 + ab[x,y], \nonumber \\ \nonumber \\
[x^2, y] &= x^2 y + y x^2 \nonumber \\
         &= x^2 y + xyx + xyx + y x^2 \nonumber \\
         &= [x,[x,y]], \mbox{ and } \nonumber \\ \label{eq2} \\
[x^2,y^2] &= [x,[x,y^2]] \nonumber \\
          &= [x,[y^2,x]] \nonumber \\
          &= [x, [y,[y,x]]]. \nonumber
\end{align}

Next, we define $C = \langle x^2 \mid x \in L_1 \rangle$ to be the 
vector subspace of $End(L)$ spanned by the squares of elements in $L_1$ and we set $S = L + C$.
Then by construction $S$ is a vector space. Equation (\ref{eq2}) asserts
that $S$ is a Lie algebra and $L$ is an ideal of $S$ with $S/L$ being abelian. Equation
(\ref{eq1}) yields that $\dim(S/L) \leq \dim(L_1)$.
\medskip

We define $S_0 = L_0+C$ and $S_1 = L_1$, and show that this 
yields a $\mathbb{Z}_2$-grading of $S$. Clearly, $[S_1, S_1] = [L_1, L_1] \subseteq 
L_0 \subseteq S_0$. Furthermore, $[S_0, S_0] \subseteq L_0$ by Equation (\ref{eq2}) and
similarly, $[S_0, S_1]  \subseteq L_1 = S_1$. It remains to show that $S_0
\cap S_1 = \{0\}$. Let $x \in S_0 \cap S_1$. Then $x \in L_1$ and $x = 
u+c$ with $u \in L_0, c \in C$. If $y \in L_0$, then $[x,y] \in L_1$
and $[x,y] = [u+c,y] = [u,y] + [c,y] \in L_0$. If $y \in L_1$, then 
$[x,y] \in L_0$ and $[x,y] = [u+c,y] = [u,y] + [c,y] \in L_1$. Since
in both cases $y \in L_0$ and $y \in L_1$, we obtain that $[x,y] = 0$,
and therefore $x \in Z(L) = \{0\}$ follows. In summary, if $x \in S_0 \cap S_1$, then $x = 0$ 
and thus $S_0 \cap S_1 = \{0\}$.
\medskip  

As $S$ has a squaring by construction, it follows that $S$ is a Lie
superalgebra. It remains to show that $S$ is simple as Lie superalgebra.
This is also shown in \cite[Th. 3.3.1]{BoLe14}. We recall a proof here 
for completeness. 
By the remarks at the end of \cite[Chapter 2]{StFa88} a minimal $2$-envelope $G$ of $L$
can be constructed as $G=L+\langle x^2 \mid x \in L \rangle$. We then have $L \subseteq S \subseteq G$. 
Let $I$ be a superideal in $S$; that is, $I$ is a Lie ideal in the Lie algebra $S$ and $s(I \cap S_1) \subseteq I$. 
Let $H$ be the ideal generated by $I$ in $G$. Then $H \cap L$ is an ideal
in $L$ and because $L$ is simple, we have $H \cap L = \{0\}$ or $H \cap L = L$.
If $H \cap L = \{0\}$, then $[H,G] \subseteq H \cap [G,G] \subseteq H \cap L = \{0\}$.
It follows that $H \subseteq Z(G) \subseteq L$, where the second inclusion follows
from Theorem 5.8 (3) in \cite[Chapter 2]{StFa88} and uses that $G$ is a minimal $2$-envelope of $L$. 
In conclusion $H=\{0\}$ and because $I\subseteq H$ also $I = \{0\}$. If $H\cap L=L$, then $L\subseteq I$ 
and thus also $C = \langle s(L_1) \rangle \subseteq \langle s(I \cap S_1) \rangle \subseteq I$ and we deduce 
$I = S$.
\end{proof}

The proof of Theorem \ref{super}
readily translates into an algorithm for constructing a {\em superization} $S$ 
based on a given $\mathbb{Z}_2$-grading of $L$; that is, an algorithm that yields an 
embedding of $L$ into a simple Lie superalgebra $S$. We determine the superizations for gradings
determined by idempotents for all known simple Lie algebras of dimension 
at most $20$ over the field $\mathbb{F}_2$ in Section \ref{results} below.


\section{Computing all subalgebras}

Let $L$ be a finite-dimensional Lie algebra over a finite field $\mathbb{F}$. Our 
aim is to introduce a practical algorithm for computing all subalgebras 
of $L$. The automorphism group $Aut(L)$ acts on the set of all subalgebras 
of $L$. Given $A \leq Aut(L)$, our method constructs orbit representatives 
under the action of $A$. We remark that in our applications $A=Aut(L)$,
but the algorithm works in the general setting.
Section \ref{results} exhibits the application of our method to the simple Lie 
algebras over $\mathbb{F}_2$ of dimension at most $16$.

\subsection{The subalgebra algorithm}

The basic idea of our method is induction. In the initial step of the 
induction, we determine the $A$-orbits of 1-dimensional subspaces of 
$L$. If $U = \langle u \rangle$ is a $1$-dimensional subspace of $L$, then
$[u,u] = 0$ and thus $U$ is a subalgebra of $L$. Hence all $1$-dimensional
subspaces are subalgebras. We denote the $A$-orbits of $1$-dimensional
subspaces with $O_1, \ldots, O_r$.

We now initialize a list $\mathcal{L}$ containing $A$-orbit representatives of all
constructed subalgebras. We choose the orbit representatives in a canonical
way by choosing the subalgebra with the lexicographically smallest 
upper triangular basis. With this convention, we obtain each $A$-orbit
of subalgebras exactly once. We only store a representative of the orbit
and we can easily check if two subalgebras belong to the same orbit. 
In the first step we add canonical $A$-orbit representatives for the
$1$-dimensional subspaces to $\mathcal{L}$.

In the induction step, we consider each subalgebra $U$ in $\mathcal{L}$ in turn.
If $\dim(U) = \dim(L)$ then $U = L$ and there is nothing to do. Thus assume 
that $\dim(U) < \dim(L)$. Then we determine $B = Stab_A(U)$ and loop over 
the orbits $O_1, \ldots, O_r$. Given $O_i$ we determine the $B$-orbits
in $O_i$. Then for each $B$-orbit representative $W$, say, we construct 
the subalgebra $V$ generated by $U$ and $W$. We add a canonical $A$-orbit
representative of $V$ to $\mathcal{L}$, if this is not already contained in $\mathcal{L}$.

\begin{remark}
We add a few remarks on efficiency of time and space of this method.
\begin{items}
\item[\rm (a)]
Storing $A$-orbit representatives instead of all subalgebras reduces the
space used to store the results.
\item[\rm (b)]
By using canonical representatives of $A$-orbits we can readily
check if a new $A$-orbit representative is already existing in $\mathcal{L}$.
\end{items}
\end{remark}

Let $U$ be a subalgebra in $\mathcal{L}$. If there exists a subalgebra $V \subseteq L$
with $U \subseteq V$, then $[V,U] \subseteq V$ and hence the quotient space $V/U$
is a submodule for the action of $ad(U)$ on $L/U$. The MeatAxe is able
to determine all submodules of a given module over a finite field and 
this is a highly efficient algorithm, provided that there are only rather
few submodules. If $L/U$ as an $ad(U)$-module has only few submodules,
then we use this alternative to construct these submodules. Once the submodules are
available, we then check which of them are subalgebras.

\subsection{Example: the $3$-dimensional simple Lie algebra}

Let $L$ be the $3$-dimensional simple Lie algebra over $\mathbb{F}_2$.
This has a basis $\{b_1, b_2, b_3\}$ with $[b_1, b_2] = b_3, [b_1, b_3]=
b_1$ and $[b_2, b_3] = b_1+b_2$. Note that this is not the standard basis
of the $3$-dimensional Zassenhaus algebra, but the basis used in the FinLie 
package. Its automorphism group $A = Aut(L)$ is a non-abelian subgroup of order $6$ 
in $GL(3,2)$. It has $3$ orbits of $1$-dimensional subspaces of $L$:
\begin{items}
\item[\rm ($O_1$)]
$\{ \langle (1,0,1) \rangle \}$,
\item[\rm ($O_2$)]
$\{ \langle (0,0,1) \rangle, \langle (1,1,0) \rangle, 
    \langle (0,1,0) \rangle \}$,
\item[\rm ($O_3$)]
$\{ \langle (0,1,1) \rangle, \langle (1,1,1) \rangle, 
    \langle (1,0,0) \rangle \}$.
\end{items}
Hence there are $7$ subalgebras of dimension $1$ falling into $3$ orbits
under $A$. The list $\mathcal{L}$ is then initialized with the representatives
$\langle (1,0,1) \rangle, \langle (0,0,1) \rangle, \langle (0,1,1) 
\rangle$. 

In the next step of the algorithm, iterated extensions are determined. This
starts with the subalgebra $U = \langle (1,0,1) \rangle$. Its stabilizer
$B$ is equal to $A$ and hence we extend $U$ twice: first with 
$W_1 = \langle (0,0,1) \rangle$ and second with $W_2 = \langle (0,1,1)
\rangle$. This yields $V_1 = \langle U, W_1 \rangle$ and $V_2 = \langle
U, W_2 \rangle$. Both are subalgebras and both have the same canonical
representative $V = \langle (1,0,0), (0,0,1) \rangle$ under the action 
of $A$. Hence they both represent the same $A$-orbit of subalgebras of
$L$ and we add its canonical representative $V$ to $\mathcal{L}$. 

In the second iteration the subalgebra $\langle (0,0,1) \rangle$ is 
extended by subalgebras in $O_2$ and $O_3$ and in the third step the
subalgebra $\langle (0,1,1) \rangle$ is extended by subalgebras in 
$O_3$. None of these extensions yields a new $A$-orbit of subalgebras,
and hence we conclude that the list $\mathcal{L}$ is already complete. 

We summarize the subalgebras of $L$ in the graph exhibited in Figure 
\ref{fig:subalgebras}. 
The vertices of the graph correspond to the subalgebras of $L$. The top vertex 
corresponds to $L$, the next layer contains the three subalgebras of dimension $2$, 
the third layer contains the seven subalgebras of dimension $1$, and at the bottom there is
the trivial subalgebra. Two subalgebras $U$ and $V$ are joined by an edge
if $U < V$ and there is no intermediate subalgebra between them. The circles 
around subalgebras of dimension $1$ and $2$ indicate $Aut(L)$-orbits.

\begin{figure}[h!] 
\begin{center}
\includegraphics[scale=0.4]{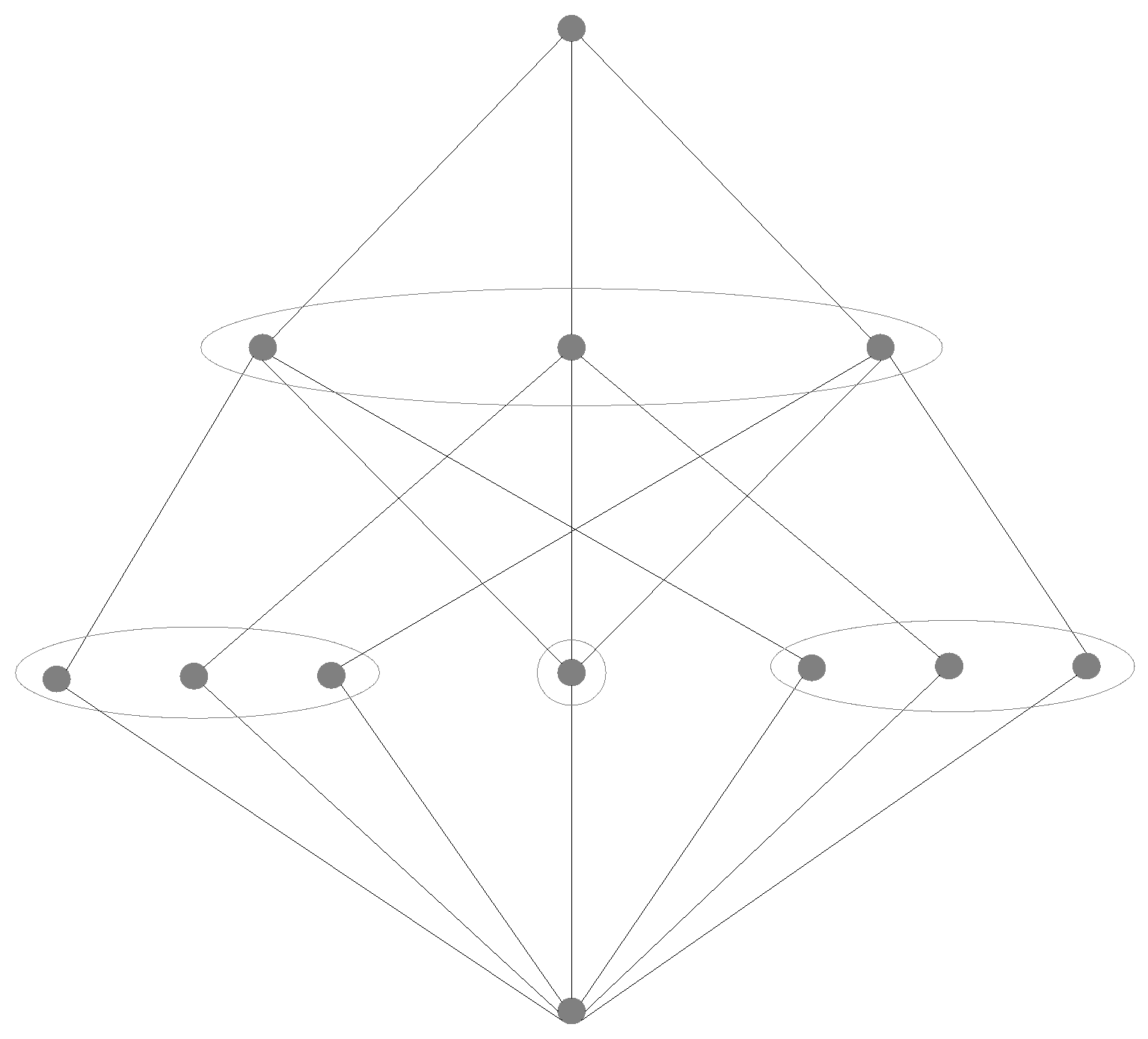}
\end{center}
\caption{The Hasse diagram of the subalgebra lattice of $L_{3,1}$.}
\label{fig:subalgebras}
\end{figure}

We note that the maximal nilpotent subalgebras of $L$ are the subalgebras
of dimension $1$. Hence $L$ is an example of a Lie algebra where the 
Cartan subalgebras are not all in one orbit under $Aut(L)$.
\subsection{Alternative approaches}

A first naive approach to compute all subalgebras is to determine all 
subspaces of $L$ as an $\mathbb{F}$-vector space and then to select those subspaces that
are closed under the Lie multiplication. This naive approach is practical 
if $|L| \leq 500$.

If $L$ has a non-trivial ideal $I$, then this can be used for an 
inductive approach. First, construct all subalgebras of the quotient 
$L/I$. Second, consider each determined subalgebra $U/I$ in turn and
determine all proper supplements to $I$ in $U$ that are subalgebras.
This idea is likely to more efficient than the method described above,
but it does not apply to the case of simple Lie algebras and this is
our desired application.


\section{Results}
\label{results}

Let $L_{d,i}$ denote the $i$-th simple Lie algebra over $\mathbb{F}_2$ of dimension 
$d$ as contained in the FinLie package. When possible, we also exhibit 
name(s) for these Lie algebras as follows:
\begin{items}
\item[$\bullet$] 
$A, B, C, D, E, F, G$ describe the simple constituents of the classical 
Lie algebras;
\item[$\bullet$] 
$W, S, H, K$ describe the simple constituents of the Lie algebras of 
Cartan type; 
\item[$\bullet$] 
$P$ describes the Hamiltonian type Lie algebras, see \cite{Lin93};
\item[$\bullet$] 
$Q$ describes the Contact type Lie algebras, see \cite{ZhLin92};
\item[$\bullet$] 
$Kap_i$ $(1 \leq i \leq 4)$ describes four series of Lie algebras constructed by Kaplansky 
\cite{Kap82};
\item[$\bullet$] 
$\Bro_i$ $(1 \leq i \leq 3)$ describes three series of Lie algebras constructed by Brown
\cite{Br95};
\item[$\bullet$] 
$V_7,V_8$ and $V_9$ are simple Lie algebras determined by Vaughan-Lee 
\cite{VL06}.
\end{items}

\subsection{Idempotents, gradings and superizations}

The following tables contain for each known simple Lie algebra $L$ 
over $\mathbb{F}_2$ up to dimension $20$ its name(s) as far as available, the 
number of $Aut(L)$-orbits of idempotents of $L$, and the dimensions of
the associated superizations. The orbits of idempotents are described
by $m \times [d_0, d_1]$, which means that there are $m$ orbits of idempotents
such that for the associated $\mathbb{Z}_2$-grading $\dim(L_0) = d_0$ and 
$\dim(L_1) = d_1$ holds. 

\begin{center}
\resizebox{\textwidth}{!}{
\begin{tabular}{|l||c|c|c|}\hline \
Lie alg & names & gradings & dim of superizations \\ 
\hline 
\hline 
$L_{3,1}$ & $W(2)$ & $1\times [1,2]$ & $5$ \\ 
\hline \hline
$L_{6,1}$ & $W(2) \otimes \mathbb{F}_4$ & $1\times [2,4]$ & $10$  \\ 
\hline \hline
$L_{7,1}$ & $W(3)$ & $2\times [3,4]$ & $9,9$ \\ 
\hline 
$L_{7,2}$ & $V_7, P(1,2)$ & $1\times [3,4]$ & $10$ \\ 
\hline \hline
$L_{8,1}$ & $A_2, W(1,1), Q(1,1,1)$ & $1\times [4,4]$ & $8$ \\ 
\hline 
$L_{8,2}$ & $V_8$ &  $1\times [4,4]$ & $8$ \\ 
\hline \hline
$L_{9,1}$ & $W(2) \otimes \mathbb{F}_8, V_9$ & $1\times [3,6]$ & $15$ \\ 
\hline \hline 
$L_{10,1}$ & $Kap_3(5)$ & $1\times [4,6]$, $1\times [6,4]$ & $12, 14$ \\ 
\hline \hline
$L_{12,1}$ & $W(2) \otimes \mathbb{F}_{16}$ & $1\times [4,8]$ & $20$ \\ 
\hline \hline
$L_{14,1}$ & $W(3) \otimes \mathbb{F}_4$ & $2\times [6,8]$ & $18, 18$ \\ 
\hline 
$L_{14,2}$ & $V_7 \otimes \mathbb{F}_4$ & $1\times [6,8]$ & $20$ \\ 
\hline 
$L_{14,3}$ & $S(2,2)$ & $3\times [6,8]$ & $16, 17, 17$ \\ 
\hline 
$L_{14,4}$ & $P(1,1,1,1), Kap_1(4)$ & $2\times [6,8]$ & $16, 17$ \\ 
\hline 
$L_{14,5}$ & {\tiny $A_3, B_3, C_3, G_2, S(1,1,1), H(1,1,1,1)$} & $1\times [6,8]$ & $14$ \\ 
\hline 
$L_{14,6}$ & $\Bro_2(1,1)$ & $1\times [6,8]$ & $16$ \\ 
\hline \hline
$L_{15,1}$ & $W(2) \otimes \mathbb{F}_{32}$ & $1\times [5,10]$ & $25$ \\ 
\hline 
$L_{15,2}$ & $W(4)$ & $4\times [7,8]$ & $17, 17, 17, 17$ \\ 
\hline 
$L_{15,3}$ & $Kap_3(6), Kap_2(4)$ & $2\times [7,8]$ & $17, 19$ \\ 
\hline 
$L_{15,4}$ & $P(2,1,1)$ & $6\times [7,8]$ & $18, 19, 19, 19, 19$ \\ 
\hline 
$L_{15,5}$ & $P(3,1)$ & $3\times [7,8]$ & $18, 18, 18$ \\ 
\hline 
$L_{15,6}$ & $P(2,2)$ & $5\times [7,8]$ & $18, 19, 19, 19, 19$ \\ 
\hline 
$L_{15,7}$ & from \cite{Eic10} & $6\times [7,8]$ & $18, 18, 19, 19, 19, 19$ \\ 
\hline 
$L_{15,8}$ & from \cite{Eic10} & $6\times [7,8]$ & $17, 17, 18, 18, 19, 19$ \\ 
\hline 
$L_{15,9}$ & new & $3\times [7,8]$ & $17, 17, 19$ \\ 
\hline 
$L_{15,10}$ & new & $8\times [7,8]$ & $17, 17, 18, 18, 19, 19, 19, 19$ \\ 
\hline 
$L_{15,11}$ & new & $3\times [7,8]$ & $17, 17, 19$ \\ 
\hline 
$L_{15,12}$ & from \cite{GriZu17},\cite{Skr98} & $4\times [7,8]$ & $18, 19, 19, 19$ \\ 
\hline \hline
$L_{16,1}$ & $W(1,1) \otimes \mathbb{F}_4, A_2 \otimes \mathbb{F}_4, V_8 \otimes \mathbb{F}_4$ & $1\times [8,8]$ & $16$ \\ 
\hline 
$L_{16,2}$ & $W(2,1), Q(2,1,1)$ & $4\times [8,8]$ & $16, 17, 17, 17$ \\ 
\hline 
$L_{16,3}$ & from \cite{Eic10} & $4\times [8,8]$ & $16, 17, 17, 17$ \\
\hline 
$L_{16,4}$ & from \cite{Eic10} & $10\times [8,8]$ & $16, 16, 17, 17, 17, 17, 17, 17, 17, 17$ \\
\hline 
$L_{16,5}$ & from \cite{Eic10} & $6\times [8,8]$ & $17, 17, 17, 17, 17, 17$ \\
\hline 
$L_{16,6}$ & from \cite{Eic10} & $4\times [8,8]$ & $16, 17, 17, 17$ \\ 
\hline \hline
$L_{18,1}$ & $W(2) \otimes \mathbb{F}_{64}$ & $1\times [6,12]$ & $30$ \\ 
\hline \hline
$L_{20,1}$ & $Kap_3(5) \otimes \mathbb{F}_4$ & $1\times [8,12]$, $1\times [12,8]$ & $24, 28$ \\ 
\hline 
\end{tabular}
}
\captionof{table}{Idempotents and superizations}
\end{center}
Note that the FinLie package contains simple Lie algebras up to dimension 
$30$ over $\mathbb{F}_2$. It is possible to compute idempotents and associated superizations 
for these algebras, but the automorphism group computation is not feasible in all cases.

\medskip
Furthermore, we remark that in all cases except for $Kap_3(5)$ (and its tensor product with
$\mathbb{F}_4$) the simple Lie algebra $L$ uniquely determines the dimensions
of $L_0$ and $L_1$ in the $\mathbb{Z}_2$-gradings associated to idempotents of $L$.

\begin{remark}
\label{tensor}
Let $L$ be a simple Lie algebra over $\mathbb{F}_2$ and let $\mathbb{E}$ be a finite field 
extension of $\mathbb{F}_2$. If $l$ is an idempotent of $L$,
then $l \otimes 1$ is an idempotent of $L \otimes_{\mathbb{F}_2} \mathbb{E}$.
\end{remark}

\subsection{Subalgebras and maximal subalgebras}

The following tables contain for each simple Lie algebra $L$ over $\mathbb{F}_2$ 
up to dimension $16$ (with the exception of $W(4) = L_{15,2}$) the number of 
$Aut(L)$-orbits of subalgebras of given dimension and the number of 
$Aut(L)$-orbits of maximal subalgebras of given dimension.

The Lie algebra $L_{15,2}$ is $W(4)$, and it is excluded in the list below.
It has more than two million orbits of subalgebras and our algorithm did 
not succeed in listing all of them.

\begin{center}
{\footnotesize \begin{tabular}{|l||r|r|r|r|r|r|r|r|r|r} 
\hline  
Lie alg &  dim 1 & dim 2 & dim 3 & dim 4 & dim 5 & dim 6 & dim 7 & dim 8 & \\ 
\hline \hline 
$L_{3,1}$ &  3 & 1 & 0 & 0 & 0 & 0 & 0 & 0 & all \\ 
\hline 
    &  0 & 1 & 0 & 0 & 0 & 0 & 0 & 0 & max \\ 
\hline 
\hline 
$L_{6,1}$ &  5 & 5 & 3 & 1 & 0 & 0 & 0 & 0 & all \\ 
\hline 
    &  0 & 0 & 1 & 1 & 0 & 0 & 0 & 0 & max \\ 
\hline 
\hline 
$L_{7,1}$ &  39 & 85 & 79 & 48 & 9 & 1 & 0 & 0 & all \\ 
\hline 
    &  0 & 0 & 1 & 1 & 0 & 1 & 0 & 0 & max \\ 
\hline 
\hline 
$L_{7,2}$ &  43 & 43 & 26 & 9 & 2 & 0 & 0 & 0 & all \\ 
\hline 
    &  0 & 0 & 4 & 5 & 2 & 0 & 0 & 0 & max \\ 
\hline 
\hline 
$L_{8,1}$ &  6 & 10 & 10 & 7 & 4 & 1 & 0 & 0 & all \\ 
\hline 
    &  0 & 0 & 0 & 0 & 1 & 1 & 0 & 0 & max \\ 
\hline 
\hline 
$L_{8,2}$ &  6 & 8 & 6 & 4 & 2 & 0 & 0 & 0 & all \\ 
\hline 
    &  0 & 0 & 0 & 1 & 2 & 0 & 0 & 0 & max \\ 
\hline 
\hline 
$L_{9,1}$ &  7 & 9 & 6 & 3 & 3 & 1 & 0 & 0 & all \\ 
\hline 
    &  0 & 0 & 1 & 0 & 0 & 1 & 0 & 0 & max \\ 
\hline 
\hline 
$L_{10,1}$ &  16 & 31 & 41 & 26 & 12 & 6 & 2 & 0 & all \\ 
\hline 
    &  0 & 0 & 0 & 0 & 0 & 2 & 2 & 0 & max \\ 
\hline 
\hline 
$L_{12,1}$ &  11 & 28 & 19 & 7 & 5 & 12 & 5 & 1 & all \\ 
\hline 
    &  0 & 0 & 0 & 0 & 0 & 1 & 0 & 1 & max \\ 
\hline 
\end{tabular}}
\captionof{table}{Subalgebras of simple Lie algebras of dimension at most $12$}
\end{center}

\begin{center}
\resizebox{\textwidth}{!}{%
\begin{tabular}{|l||r|r|r|r|r|r|r|r|r|r|r|r|r|r} 
\hline 
Lie alg &  dim 1 & dim 2 & dim 3 & dim 4 & dim 5 & dim 6 & dim 7 & dim 8 & dim 9 & dim 10 & dim 11 & dim 12 & \\ 
\hline \hline 
$L_{14,1}$ &  211 & 2712 & 8011 & 9548 & 9827 & 8345 & 6564 & 2778 & 316 & 28 & 5 & 1 & all \\ 
\hline 
    &  0 & 0 & 0 & 0 & 0 & 1 & 1 & 1 & 0 & 0 & 0 & 1 & max \\ 
\hline 
\hline 
$L_{14,2}$ &  237 & 481 & 532 & 251 & 116 & 58 & 36 & 19 & 4 & 2 & 0 & 0 & all \\ 
\hline 
    &  0 & 0 & 0 & 0 & 0 & 7 & 1 & 6 & 0 & 2 & 0 & 0 & max \\ 
\hline 
\hline 
$L_{14,3}$ &  135 & 790 & 1988 & 2545 & 2315 & 1489 & 822 & 298 & 56 & 11 & 2 & 1 & all \\ 
\hline 
    &  0 & 0 & 0 & 0 & 0 & 1 & 3 & 1 & 0 & 0 & 1 & 1 & max \\ 
\hline 
\hline 
$L_{14,4}$ &  78 & 289 & 538 & 545 & 360 & 204 & 139 & 56 & 14 & 3 & 1 & 0 & all \\ 
\hline 
    &  0 & 0 & 0 & 0 & 0 & 1 & 3 & 2 & 0 & 3 & 1 & 0 & max \\ 
\hline 
\hline 
$L_{14,5}$ &  6 & 14 & 25 & 32 & 28 & 19 & 13 & 10 & 4 & 1 & 1 & 0 & all \\ 
\hline 
    &  0 & 0 & 0 & 0 & 0 & 0 & 0 & 1 & 0 & 1 & 1 & 0 & max \\ 
\hline 
\hline 
$L_{14,6}$ &  19 & 70 & 143 & 171 & 126 & 81 & 53 & 30 & 10 & 1 & 2 & 0 & all \\ 
\hline 
    &  0 & 0 & 0 & 0 & 0 & 1 & 1 & 0 & 0 & 1 & 2 & 0 & max \\ 
\hline 
\end{tabular}
}
\captionof{table}{Subalgebras of simple Lie algebras of dimension $14$}
\end{center}

\begin{center}
\resizebox{\textwidth}{!}{%
\begin{tabular}{|l||r|r|r|r|r|r|r|r|r|r|r|r|r|r|r} 
\hline 
Lie alg &  dim 1 & dim 2 & dim 3 & dim 4 & dim 5 & dim 6 & dim 7 & dim 8 & dim 9 & dim 10 & dim 11 & dim 12 & dim 13 & \\ 
\hline \hline 
$L_{15,1}$ &  15 & 67 & 72 & 19 & 5 & 7 & 31 & 31 & 7 & 1 & 0 & 0 & 0 & all \\ 
\hline 
    &  0 & 0 & 1 & 0 & 0 & 0 & 0 & 0 & 0 & 1 & 0 & 0 & 0 & max \\ 
\hline 
\hline 
$L_{15,3}$ &  37 & 156 & 318 & 447 & 366 & 257 & 159 & 94 & 30 & 10 & 3 & 1 & 0 & all \\ 
\hline 
    &  0 & 0 & 0 & 0 & 0 & 2 & 0 & 0 & 0 & 1 & 2 & 1 & 0 & max \\ 
\hline 
\hline 
$L_{15,4}$ &  911 & 3177 & 5925 & 7027 & 5250 & 3260 & 1714 & 721 & 167 & 34 & 10 & 2 & 0 & all \\ 
\hline 
    &  0 & 0 & 14 & 0 & 0 & 1 & 15 & 13 & 2 & 0 & 4 & 2 & 0 & max \\ 
\hline 
\hline 
$L_{15,5}$ &  455 & 4199 & 14178 & 23832 & 26523 & 22453 & 13857 & 5663 & 1555 & 316 & 53 & 8 & 2 & all \\ 
\hline 
    &  0 & 0 & 2 & 0 & 0 & 0 & 3 & 4 & 0 & 0 & 0 & 0 & 2 & max \\ 
\hline 
\hline 
$L_{15,6}$ &  511 & 2975 & 7644 & 11384 & 9992 & 6933 & 4018 & 1948 & 625 & 143 & 24 & 5 & 1 & all \\ 
\hline 
    &  0 & 0 & 5 & 0 & 0 & 2 & 5 & 4 & 0 & 0 & 1 & 1 & 1 & max \\ 
\hline 
\hline 
$L_{15,7}$ &  1663 & 4823 & 7807 & 8280 & 5229 & 2978 & 1547 & 826 & 162 & 21 & 4 & 1 & 0 & all \\ 
\hline 
    &  0 & 0 & 26 & 0 & 0 & 1 & 30 & 22 & 7 & 1 & 3 & 1 & 0 & max \\ 
\hline 
\hline 
$L_{15,8}$ &  475 & 1885 & 3410 & 4010 & 2720 & 1631 & 861 & 472 & 97 & 18 & 5 & 1 & 0 & all \\ 
\hline 
    &  0 & 0 & 5 & 0 & 0 & 3 & 6 & 6 & 2 & 1 & 4 & 1 & 0 & max \\ 
\hline 
\hline 
$L_{15,9}$ &  117 & 542 & 1146 & 1603 & 1245 & 807 & 464 & 273 & 81 & 17 & 4 & 2 & 0 & all \\ 
\hline 
    &  0 & 0 & 1 & 0 & 0 & 0 & 2 & 1 & 0 & 0 & 2 & 2 & 0 & max \\ 
\hline 
\hline 
$L_{15,10}$ &  491 & 2307 & 4946 & 6504 & 5070 & 3214 & 1666 & 706 & 152 & 34 & 11 & 2 & 0 & all \\ 
\hline 
    &  0 & 0 & 6 & 0 & 0 & 1 & 9 & 6 & 1 & 1 & 3 & 2 & 0 & max \\ 
\hline 
\hline 
$L_{15,11}$ &  77 & 304 & 608 & 781 & 575 & 385 & 223 & 125 & 41 & 10 & 2 & 1 & 0 & all \\ 
\hline 
    &  0 & 0 & 1 & 0 & 0 & 1 & 3 & 1 & 0 & 0 & 1 & 1 & 0 & max \\ 
\hline 
\hline 
$L_{15,12}$ &  687 & 3449 & 8572 & 11602 & 8502 & 4533 & 2120 & 895 & 206 & 30 & 9 & 3 & 0 & all \\ 
\hline 
    &  0 & 0 & 6 & 0 & 0 & 0 & 10 & 8 & 1 & 0 & 2 & 3 & 0 & max \\ 
\hline 
\end{tabular}
}
\captionof{table}{Subalgebras of simple Lie algebras of dimension $15$ except $W(4) = L_{15,2}$}
\end{center}

\begin{center}
\resizebox{\textwidth}{!}{%
\begin{tabular}{|l||r|r|r|r|r|r|r|r|r|r|r|r|r|r|r|r} 
\hline 
Lie alg &  dim 1 & dim 2 & dim 3 & dim 4 & dim 5 & dim 6 & dim 7 & dim 8 & dim 9 & dim 10 & dim 11 & dim 12 & dim 13 & dim 14 & \\ 
\hline \hline 
$L_{16,1}$ &  13 & 51 & 56 & 70 & 49 & 52 & 41 & 42 & 16 & 6 & 2 & 1 & 0 & 0 & all \\ 
\hline 
    &  0 & 0 & 0 & 0 & 0 & 0 & 0 & 2 & 0 & 1 & 0 & 1 & 0 & 0 & max \\ 
\hline 
\hline 
$L_{16,2}$ &  157 & 1445 & 5214 & 10302 & 11518 & 10008 & 6604 & 3878 & 1414 & 351 & 68 & 19 & 5 & 1 & all \\ 
\hline 
    &  0 & 0 & 0 & 0 & 0 & 0 & 0 & 2 & 4 & 1 & 0 & 1 & 0 & 1 & max \\ 
\hline 
\hline 
$L_{16,3}$ &  168 & 1086 & 2658 & 4248 & 3838 & 2485 & 1420 & 879 & 344 & 74 & 16 & 5 & 2 & 0 & all \\ 
\hline 
    &  0 & 0 & 0 & 0 & 0 & 0 & 0 & 2 & 3 & 1 & 0 & 3 & 2 & 0 & max \\ 
\hline 
\hline 
$L_{16,4}$ &  495 & 2988 & 6900 & 10355 & 8700 & 5457 & 2968 & 1752 & 651 & 109 & 18 & 6 & 1 & 0 & all \\ 
\hline 
    &  0 & 0 & 0 & 1 & 0 & 0 & 0 & 7 & 9 & 3 & 0 & 4 & 1 & 0 & max \\ 
\hline 
\hline 
$L_{16,5}$ &  379 & 2267 & 5857 & 9320 & 8037 & 5220 & 2930 & 1832 & 696 & 137 & 23 & 8 & 2 & 0 & all \\ 
\hline 
    &  0 & 0 & 0 & 2 & 0 & 0 & 0 & 7 & 8 & 1 & 0 & 2 & 2 & 0 & max \\ 
\hline 
\hline 
$L_{16,6}$ &  297 & 972 & 1311 & 1397 & 791 & 232 & 89 & 86 & 34 & 5 & 1 & 1 & 0 & 0 & all \\ 
\hline 
    &  0 & 0 & 0 & 3 & 0 & 0 & 0 & 9 & 6 & 4 & 1 & 1 & 0 & 0 & max \\ 
\hline 
\end{tabular}
}
\captionof{table}{Subalgebras of simple Lie algebras of dimension $16$}
\end{center}


\subsection{Subquotients}

The following tables exhibit for each simple Lie algebra $L$ over $\mathbb{F}_2$ 
up to dimension $16$ (with the exception of $W(4) = L_{15,2}$) the simple 
Lie algebras which arise as proper subquotients. An entry ``+'' in the
row of $L_{d,i}$ and column $L_{f,j}$ indicates that $L_{f,j}$ embeds as
a proper subquotient into $L_{d,i}$.

\begin{center}

\begin{center}
{\footnotesize
\begin{tabular}{|c||c|c|c|c|c|c|c|c|c|} \hline 
& $L_{3,1}$ & $L_{6,1}$ & $L_{7,1}$ & $L_{7,2}$ & $L_{8,1}$ & $L_{8,2}$ & $L_{9,1}$ & $L_{10,1}$ \\ \hline \hline
$L_{3,1}$ &   & & & & & & & \\ \hline
$L_{6,1}$ & + & & & & & & & \\ \hline
$L_{7,1}$ & + & & & & & & & \\ \hline
$L_{7,2}$ & + & & & & & & & \\ \hline
$L_{8,1}$ & + & & & & & & & \\ \hline
$L_{8,2}$ & + & & & & & & & \\ \hline
$L_{9,1}$ & + & & & & & & & \\ \hline
$L_{10,1}$ & + & & & & & & & \\ \hline
$L_{12,1}$ & + & $+$ & & & & & & \\ \hline
$L_{14,1}$ & + & $+$ & $+$ & & & & & \\ \hline
$L_{14,2}$ & + & $+$ & & $+$ & & & & \\ \hline
$L_{14,3}$ & + & & $+$ & $+$ & & & & \\ \hline
$L_{14,4}$ & + & & $+$ & $+$ & $+$ & & & \\ \hline
$L_{14,5}$ & + & & & & $+$ & $+$ & & \\ \hline
$L_{14,6}$ & + & & $+$ & & $+$ & & & \\ \hline
$L_{15,1}$ & + & & & & & & & \\ \hline
$L_{15,3}$ & + & $+$ & & & $+$ & & & $+$ \\ \hline
$L_{15,4}$ & + & & $+$ & $+$ & & & & \\ \hline
$L_{15,5}$ & + & & $+$ & $+$ & & & & \\ \hline
$L_{15,6}$ & + & $+$ & $+$ & $+$ & & & & \\ \hline
$L_{15,7}$ & + & $+$ & $+$ & $+$ & $+$ & & & \\ \hline
$L_{15,8}$ & + & $+$ & $+$ & $+$ & $+$ & & & $+$ \\ \hline
$L_{15,9}$ & + & & $+$ & & $+$ & & & \\ \hline
$L_{15,10}$ & + & & $+$ & $+$ & & & & $+$ \\ \hline
$L_{15,11}$ & + & $+$ & $+$ & & $+$ & & & \\ \hline
$L_{15,12}$ & + & & $+$ & $+$ & & & & \\ \hline
$L_{16,1}$ & + & $+$ & & & $+$ & $+$ & & \\ \hline
$L_{16,2}$ & + & & $+$ & $+$ & $+$ & & & \\ \hline
$L_{16,3}$ & + & & $+$ & $+$ & $+$ & & & \\ \hline
$L_{16,4}$ & + & & $+$ & $+$ & $+$ & $+$ & & \\ \hline
$L_{16,5}$ & + & & $+$ & $+$ & $+$ & $+$ & & \\ \hline
$L_{16,6}$ & + & & $+$ & $+$ & $+$ & $+$ & & \\ \hline
\end{tabular}}
\end{center}

\captionof{table}{Simple subquotients of simple Lie algebras}
\end{center}

\newpage
\section{Appendix}\label{appendix}
We list explicit generators for the three new simple Lie algebras that 
have been found using the random search methods described in \cite{Eic10}. 
We note that the Lie algebras described in \cite{GriZu17} are 
isomorphic to the Lie algebra obtained in \cite{Skr98} over $\mathbb{F}_2$ and
hence are already contained in our list of known simple Lie algebras.
\medskip

\setcounter{MaxMatrixCols}{15}

Explicit generators for the Lie algebra number $9$ of dimension $15$:
\begin{center}
{\tiny
\setlength\arraycolsep{2pt}
$\begin{pmatrix} 
\cdot & \cdot & \cdot & 1 & \cdot & 1 & \cdot & 1 & \cdot & \cdot & \cdot & 1 & \cdot & \cdot & 1 \\ 
\cdot & \cdot & 1 & 1 & 1 & \cdot & 1 & 1 & \cdot & \cdot & 1 & \cdot & 1 & 1 & 1 \\ 
\cdot & 1 & 1 & \cdot & 1 & 1 & 1 & \cdot & \cdot & 1 & \cdot & \cdot & \cdot & 1 & 1 \\ 
\cdot & \cdot & 1 & \cdot & \cdot & 1 & \cdot & 1 & \cdot & \cdot & 1 & \cdot & \cdot & \cdot & \cdot \\ 
\cdot & \cdot & \cdot & 1 & \cdot & \cdot & \cdot & 1 & \cdot & \cdot & \cdot & 1 & \cdot & 1 & \cdot \\ 
\cdot & \cdot & \cdot & 1 & 1 & 1 & 1 & 1 & \cdot & \cdot & \cdot & 1 & \cdot & \cdot & \cdot \\ 
\cdot & \cdot & \cdot & \cdot & 1 & \cdot & 1 & \cdot & \cdot & 1 & 1 & 1 & \cdot & \cdot & \cdot \\ 
\cdot & \cdot & \cdot & \cdot & \cdot & 1 & 1 & \cdot & 1 & \cdot & \cdot & 1 & 1 & \cdot & \cdot \\ 
\cdot & \cdot & \cdot & \cdot & \cdot & \cdot & 1 & 1 & \cdot & 1 & 1 & 1 & \cdot & 1 & 1 \\ 
\cdot & \cdot & \cdot & \cdot & \cdot & \cdot & 1 & \cdot & 1 & \cdot & \cdot & \cdot & 1 & \cdot & \cdot \\ 
\cdot & \cdot & \cdot & \cdot & \cdot & \cdot & 1 & 1 & 1 & \cdot & 1 & \cdot & \cdot & 1 & \cdot \\ 
\cdot & \cdot & \cdot & \cdot & \cdot & 1 & 1 & \cdot & \cdot & \cdot & 1 & \cdot & 1 & \cdot & \cdot \\ 
\cdot & \cdot & \cdot & \cdot & \cdot & \cdot & 1 & 1 & 1 & \cdot & 1 & 1 & 1 & 1 & 1 \\ 
\cdot & \cdot & \cdot & \cdot & \cdot & \cdot & \cdot & 1 & \cdot & \cdot & \cdot & \cdot & 1 & 1 & 1 \\ 
\cdot & \cdot & \cdot & \cdot & \cdot & \cdot & \cdot & \cdot & \cdot & \cdot & 1 & 1 & 1 & 1 & \cdot \end{pmatrix}$,
$\begin{pmatrix}
\cdot & \cdot & \cdot & \cdot & \cdot & \cdot & 1 & \cdot & 1 & \cdot & 1 & \cdot & \cdot & 1 & 1 \\ 
\cdot & \cdot & 1 & \cdot & \cdot & 1 & \cdot & \cdot & 1 & 1 & \cdot & \cdot & \cdot & \cdot & 1 \\ 
1 & \cdot & \cdot & \cdot & 1 & \cdot & 1 & 1 & 1 & 1 & 1 & 1 & \cdot & \cdot & 1 \\ 
\cdot & \cdot & \cdot & \cdot & \cdot & 1 & 1 & \cdot & \cdot & 1 & 1 & \cdot & 1 & 1 & 1 \\ 
\cdot & \cdot & 1 & 1 & \cdot & \cdot & 1 & \cdot & 1 & 1 & \cdot & 1 & 1 & 1 & 1 \\ 
\cdot & \cdot & \cdot & 1 & 1 & \cdot & \cdot & 1 & 1 & 1 & 1 & 1 & 1 & \cdot & 1 \\ 
\cdot & \cdot & \cdot & 1 & 1 & \cdot & \cdot & 1 & \cdot & 1 & 1 & \cdot & \cdot & 1 & \cdot \\ 
\cdot & \cdot & \cdot & \cdot & 1 & \cdot & \cdot & \cdot & 1 & 1 & 1 & \cdot & \cdot & \cdot & \cdot \\ 
\cdot & \cdot & \cdot & \cdot & \cdot & \cdot & \cdot & \cdot & \cdot & 1 & 1 & 1 & \cdot & \cdot & 1 \\ 
\cdot & \cdot & \cdot & \cdot & \cdot & 1 & 1 & \cdot & \cdot & 1 & 1 & 1 & \cdot & 1 & \cdot \\ 
\cdot & \cdot & \cdot & \cdot & \cdot & 1 & 1 & \cdot & \cdot & 1 & \cdot & 1 & 1 & \cdot & 1 \\ 
\cdot & \cdot & \cdot & \cdot & \cdot & \cdot & \cdot & 1 & 1 & 1 & \cdot & \cdot & \cdot & \cdot & \cdot \\ 
\cdot & \cdot & \cdot & \cdot & \cdot & 1 & 1 & \cdot & \cdot & 1 & 1 & \cdot & 1 & \cdot & \cdot \\ 
\cdot & \cdot & \cdot & \cdot & \cdot & \cdot & 1 & \cdot & 1 & 1 & 1 & \cdot & 1 & \cdot & 1 \\ 
\cdot & \cdot & \cdot & \cdot & \cdot & \cdot & \cdot & 1 & 1 & \cdot & 1 & 1 & 1 & \cdot & \cdot \end{pmatrix}$}
\end{center}

Explicit generators for the Lie algebra number $10$ of dimension $15$:
\begin{center}
{\tiny
\setlength\arraycolsep{2pt}
$\begin{pmatrix}\cdot & \cdot & \cdot & \cdot & \cdot & \cdot & \cdot & 1 & \cdot & 1 & 1 & \cdot & 1 & 1 & 1 \\ 
\cdot & \cdot & 1 & \cdot & \cdot & \cdot & 1 & \cdot & \cdot & \cdot & 1 & \cdot & 1 & 1 & \cdot \\ 
\cdot & 1 & \cdot & 1 & 1 & 1 & \cdot & 1 & \cdot & 1 & \cdot & \cdot & \cdot & 1 & 1 \\ 
\cdot & \cdot & 1 & \cdot & 1 & \cdot & \cdot & \cdot & 1 & 1 & \cdot & 1 & \cdot & 1 & \cdot \\ 
\cdot & \cdot & \cdot & \cdot & 1 & \cdot & \cdot & \cdot & \cdot & \cdot & 1 & 1 & 1 & 1 & \cdot \\ 
\cdot & \cdot & \cdot & 1 & \cdot & \cdot & \cdot & \cdot & \cdot & \cdot & 1 & 1 & \cdot & \cdot & \cdot \\ 
\cdot & \cdot & \cdot & \cdot & 1 & 1 & 1 & 1 & 1 & \cdot & \cdot & \cdot & 1 & \cdot & 1 \\ 
\cdot & \cdot & \cdot & \cdot & \cdot & 1 & 1 & \cdot & 1 & 1 & \cdot & \cdot & 1 & 1 & 1 \\ 
\cdot & \cdot & \cdot & \cdot & \cdot & \cdot & \cdot & \cdot & 1 & 1 & 1 & 1 & 1 & 1 & \cdot \\ 
\cdot & \cdot & \cdot & \cdot & \cdot & \cdot & \cdot & \cdot & \cdot & 1 & 1 & 1 & \cdot & \cdot & \cdot \\ 
\cdot & \cdot & \cdot & \cdot & \cdot & \cdot & 1 & \cdot & \cdot & \cdot & \cdot & \cdot & 1 & \cdot & \cdot \\ 
\cdot & \cdot & \cdot & \cdot & \cdot & 1 & 1 & 1 & 1 & 1 & 1 & \cdot & \cdot & 1 & 1 \\ 
\cdot & \cdot & \cdot & \cdot & \cdot & \cdot & 1 & 1 & \cdot & \cdot & 1 & \cdot & \cdot & \cdot & \cdot \\ 
\cdot & \cdot & \cdot & \cdot & \cdot & \cdot & \cdot & 1 & \cdot & \cdot & \cdot & 1 & 1 & \cdot & \cdot \\ 
\cdot & \cdot & \cdot & \cdot & \cdot & \cdot & \cdot & \cdot & \cdot & 1 & 1 & 1 & \cdot & 1 & \cdot \end{pmatrix}$,
$\begin{pmatrix}\cdot & \cdot & \cdot & \cdot & 1 & 1 & \cdot & \cdot & \cdot & 1 & 1 & \cdot & 1 & 1 & 1 \\ 
\cdot & \cdot & 1 & 1 & 1 & 1 & \cdot & \cdot & 1 & 1 & 1 & 1 & \cdot & 1 & 1 \\ 
1 & \cdot & 1 & 1 & 1 & \cdot & \cdot & \cdot & \cdot & 1 & \cdot & \cdot & \cdot & 1 & 1 \\ 
\cdot & \cdot & 1 & \cdot & \cdot & \cdot & \cdot & 1 & 1 & \cdot & 1 & \cdot & 1 & 1 & \cdot \\ 
\cdot & \cdot & 1 & 1 & 1 & 1 & \cdot & \cdot & \cdot & \cdot & 1 & \cdot & 1 & \cdot & \cdot \\ 
\cdot & \cdot & \cdot & 1 & 1 & 1 & \cdot & \cdot & \cdot & 1 & \cdot & 1 & \cdot & 1 & \cdot \\ 
\cdot & \cdot & \cdot & 1 & \cdot & 1 & 1 & \cdot & \cdot & \cdot & 1 & \cdot & \cdot & \cdot & \cdot \\ 
\cdot & \cdot & \cdot & \cdot & 1 & \cdot & \cdot & 1 & \cdot & 1 & 1 & 1 & 1 & \cdot & \cdot \\ 
\cdot & \cdot & \cdot & \cdot & \cdot & 1 & \cdot & \cdot & 1 & 1 & 1 & \cdot & \cdot & \cdot & \cdot \\ 
\cdot & \cdot & \cdot & \cdot & \cdot & \cdot & 1 & \cdot & 1 & 1 & 1 & \cdot & 1 & \cdot & \cdot \\ 
\cdot & \cdot & \cdot & \cdot & \cdot & 1 & \cdot & \cdot & \cdot & \cdot & 1 & \cdot & 1 & 1 & \cdot \\ 
\cdot & \cdot & \cdot & \cdot & \cdot & \cdot & \cdot & \cdot & \cdot & \cdot & \cdot & \cdot & 1 & \cdot & \cdot \\ 
\cdot & \cdot & \cdot & \cdot & \cdot & 1 & 1 & 1 & \cdot & \cdot & 1 & \cdot & 1 & \cdot & \cdot \\ 
\cdot & \cdot & \cdot & \cdot & \cdot & \cdot & 1 & 1 & \cdot & \cdot & 1 & 1 & 1 & \cdot & \cdot \\ 
\cdot & \cdot & \cdot & \cdot & \cdot & \cdot & \cdot & 1 & 1 & 1 & 1 & \cdot & 1 & 1 & 1 \end{pmatrix}$}
\end{center}

Explicit generators for the Lie algebra number $11$ of dimension $15$:
\begin{center}
{\tiny
\setlength\arraycolsep{2pt}
$\begin{pmatrix}\cdot & 1 & 1 & 1 & \cdot & 1 & 1 & 1 & \cdot & 1 & 1 & \cdot & 1 & 1 & \cdot \\ 
\cdot & 1 & \cdot & 1 & 1 & \cdot & \cdot & 1 & 1 & \cdot & 1 & \cdot & \cdot & 1 & \cdot \\ 
\cdot & 1 & 1 & 1 & \cdot & 1 & 1 & \cdot & 1 & \cdot & 1 & 1 & 1 & \cdot & \cdot \\ 
\cdot & \cdot & 1 & 1 & 1 & \cdot & 1 & 1 & \cdot & 1 & \cdot & \cdot & 1 & 1 & 1 \\ 
\cdot & \cdot & \cdot & 1 & \cdot & 1 & \cdot & 1 & 1 & \cdot & \cdot & \cdot & 1 & \cdot & \cdot \\ 
\cdot & \cdot & \cdot & 1 & 1 & 1 & \cdot & \cdot & 1 & \cdot & 1 & \cdot & \cdot & \cdot & \cdot \\ 
\cdot & \cdot & \cdot & \cdot & 1 & 1 & \cdot & \cdot & \cdot & \cdot & \cdot & \cdot & 1 & 1 & \cdot \\ 
\cdot & \cdot & \cdot & \cdot & \cdot & \cdot & \cdot & \cdot & 1 & \cdot & 1 & \cdot & \cdot & \cdot & \cdot \\ 
\cdot & \cdot & \cdot & \cdot & \cdot & 1 & \cdot & 1 & 1 & \cdot & 1 & \cdot & 1 & \cdot & 1 \\ 
\cdot & \cdot & \cdot & \cdot & \cdot & 1 & 1 & \cdot & \cdot & 1 & \cdot & 1 & 1 & 1 & \cdot \\ 
\cdot & \cdot & \cdot & \cdot & \cdot & 1 & 1 & 1 & \cdot & \cdot & 1 & 1 & 1 & \cdot & \cdot \\ 
\cdot & \cdot & \cdot & \cdot & \cdot & 1 & \cdot & 1 & \cdot & \cdot & 1 & \cdot & \cdot & \cdot & \cdot \\ 
\cdot & \cdot & \cdot & \cdot & \cdot & \cdot & 1 & 1 & \cdot & 1 & 1 & 1 & 1 & 1 & \cdot \\ 
\cdot & \cdot & \cdot & \cdot & \cdot & \cdot & \cdot & 1 & \cdot & 1 & 1 & \cdot & \cdot & 1 & \cdot \\ 
\cdot & \cdot & \cdot & \cdot & \cdot & \cdot & \cdot & \cdot & 1 & \cdot & \cdot & 1 & 1 & \cdot & 1 \end{pmatrix}$,
$\begin{pmatrix}1 & \cdot & 1 & \cdot & 1 & \cdot & \cdot & \cdot & \cdot & 1 & \cdot & 1 & \cdot & 1 & 1 \\ 
1 & \cdot & 1 & 1 & \cdot & 1 & 1 & 1 & \cdot & 1 & 1 & \cdot & \cdot & 1 & \cdot \\ 
1 & \cdot & 1 & \cdot & 1 & \cdot & 1 & 1 & \cdot & 1 & \cdot & 1 & \cdot & \cdot & \cdot \\ 
\cdot & \cdot & \cdot & 1 & 1 & 1 & \cdot & 1 & 1 & 1 & 1 & \cdot & 1 & 1 & \cdot \\ 
\cdot & \cdot & 1 & \cdot & \cdot & 1 & 1 & 1 & \cdot & \cdot & \cdot & 1 & 1 & \cdot & \cdot \\ 
\cdot & \cdot & \cdot & 1 & 1 & \cdot & 1 & 1 & 1 & 1 & \cdot & 1 & 1 & \cdot & 1 \\ 
\cdot & \cdot & \cdot & 1 & 1 & \cdot & 1 & \cdot & \cdot & \cdot & \cdot & 1 & 1 & 1 & \cdot \\ 
\cdot & \cdot & \cdot & \cdot & 1 & 1 & \cdot & 1 & \cdot & \cdot & 1 & 1 & \cdot & 1 & \cdot \\ 
\cdot & \cdot & \cdot & \cdot & \cdot & \cdot & 1 & \cdot & 1 & \cdot & \cdot & 1 & \cdot & 1 & 1 \\ 
\cdot & \cdot & \cdot & \cdot & \cdot & \cdot & 1 & \cdot & 1 & \cdot & 1 & 1 & 1 & 1 & \cdot \\ 
\cdot & \cdot & \cdot & \cdot & \cdot & \cdot & 1 & 1 & 1 & \cdot & 1 & 1 & \cdot & \cdot & 1 \\ 
\cdot & \cdot & \cdot & \cdot & \cdot & 1 & 1 & 1 & \cdot & 1 & 1 & 1 & 1 & \cdot & 1 \\ 
\cdot & \cdot & \cdot & \cdot & \cdot & 1 & 1 & 1 & \cdot & 1 & 1 & \cdot & \cdot & 1 & \cdot \\ 
\cdot & \cdot & \cdot & \cdot & \cdot & \cdot & 1 & \cdot & 1 & 1 & \cdot & \cdot & \cdot & \cdot & 1 \\ 
\cdot & \cdot & \cdot & \cdot & \cdot & \cdot & \cdot & 1 & 1 & 1 & \cdot & 1 & \cdot & 1 & \cdot \end{pmatrix}$}
\end{center}


\EditInfo{August 30, 2021}{January 6, 2022}{Friedrich Wagemann}

\end{paper}
\begin{references}

\refer{Paper}{BeBou18}
    \Rauthor{Benayadi, S.; Bouarroudj, S.}
    \Rtitle{Double extensions of {L}ie superalgebras in characteristic 2 with nondegenerate invariant supersymmetric bilinear form}
    \Rjournal{J. Algebra}
    \Rvolume{510}
    \Ryear{2018}
    \Rpages{141-179}

\refer{Arxiv}{BoLe14}
    \Rauthor{Bouarroudj, S.; Lebedev, A.; Leites, D.; Shchepochkina, I.} 
    \Rtitle{Classification of simple {L}ie superalgebras in characteristic $2$}
    \Rarxivid{ArXiv:~1407.1695}
    
\refer{Paper}{Br95}
    \Rauthor{Brown, G.}
    \Rtitle{Families of simple {L}ie algebras of characteristic two}
    \Rjournal{Comm. Algebra}
    \Rvolume {23}
    \Ryear {1995}
    \Rnumber {3}
    \Rpages {941-954}

\refer{Paper}{Cal14}
   \Rauthor{Calderon-Martin, A.}
   \Rtitle{Lie algebras with a set grading}
   \Rjournal{Lin. Alg. Appl.}
   \Rvolume{452}
   \Ryear{2014}
   \Rpages{7-20}

\refer{Paper}{Eic10}
    \Rauthor{Eick, B.}
    \Rtitle{Some new simple {L}ie algebras in characteristic 2}
    \Rjournal{J. Symbolic Comput.}
    \Rvolume{45}
    \Ryear{2010}
    \Rnumber{9}
    \Rpages{943-951}
   
\refer{Other}{finlie}
    \Rauthor{Eick, B.}
    \Rtitle{FinLie - Computing with finite Lie algebras (A (non-refereed) GAP 4 package)}
    \Ryear{2015}  

\refer{Book}{ElKo13}
    \Rtitle{Gradings on Simple Lie Algebras}
    \Rauthor{Elduque, A.; Kochetov, M.}
    \Ryear{2013}
    \Rpublisher{American Mathematical Society}

\refer{Paper}{GriZu17}
    \Rauthor{Grishkov, A.; Zusmanovich, P.}
    \Rtitle{Deformations of current {L}ie algebras. I. Small algebras in characteristic 2}
    \Rjournal{Journal of Algebra}
    \Rvolume{473}
    \Rpages{513-544}
    \Ryear{2017}

\refer{Paper}{Hi84}
    \Rauthor{Hiss, G.}
    \Rtitle{Die adjungierten {D}arstellungen der {C}hevalley-{G}ruppen}
    \Rjournal{Arch. Math. (Basel)}
    \Rvolume{42}
    \Ryear{1984}
    \Rnumber{5}
    \Rpages{408-416}

\refer{Paper}{Ho82}
    \Rauthor{Hogeweij, G. M. D.}
    \Rtitle{Almost-classical {L}ie algebras. {I}, {II}}
    \Rjournal{Nederl. Akad. Wetensch. Indag. Math.}
    \Rvolume{44}
    \Ryear{1982}
    \Rnumber{4}
    \Rpages{441-460}

\refer{Paper}{Kac77}
    \Rauthor{Kac, V. G.}
    \Rtitle{Lie superalgebras}
    \Rjournal{Advances in Mathematics}
    \Rvolume{26}
    \Rnumber{1}
    \Rpages{8-96}
    \Ryear{1977}
   
\refer{Book}{Kap82}
    \Rauthor{Kaplansky, I.}
    \Rtitle{Some simple {L}ie algebras of characteristic $2$}
    \Rbooktitle{Lie algebras and related topics ({N}ew {B}runswick, {N}.{J}., 1981)}
    \Rseries{Lecture Notes in Math.}
    \Rvolume{933}
    \Rpublisher{Springer, Berlin-New York}
    \Ryear{1982}

\refer{Paper}{KrLe18}
    \Rauthor{Krutov, A.; Lebedev, A.}
    \Rtitle{On gradings modulo 2 of simple {L}ie algebras in characteristic 2}
    \Rjournal{SIGMA Symmetry Integrability Geom. Methods Appl.}
    \Rvolume{14}
    \Rnumber{130 (27)}
    \Ryear{2018}
    
\refer{Paper}{Lin93}
    \Rauthor{Lin, L.}
    \Rtitle{Nonalternating Hamiltonian algebra $P(n, m)$ of characteristic two}
    \Rjournal{Comm. Algebra}
    \Rvolume{21}
    \Ryear{1993}
    \Rnumber{2}
    \Rpages{399-411}
    
\refer{Paper}{PZa89}
    \Rauthor{Patera, J.; Zassenhaus, H.}
    \Rtitle{On {L}ie gradings. {I}}
    \Rjournal{Linear Algebra Appl.}
    \Rvolume{112}
    \Ryear{1989}
    \Rpages{87-159}

\refer{Incollection}{PSt06}
    \Rauthor{Premet, A.; Strade, H.}
    \Rtitle{Classification of finite dimensional simple {L}ie algebras in prime characteristics}
    \Rbooktitle{Representations of algebraic groups, quantum groups, and Lie  algebras}
    \Rseries{Contemp. Math.}
    \Rvolume{413}
    \Rpages{185-214}
    \Rpublisher{Amer. Math. Soc.}
    \Ryear{2006}

\refer{Paper}{Skr98}
    \Rauthor{Skryabin, S.}
    \Rtitle{Toral rank one simple {L}ie algebras of low characteristics}
    \Rjournal{J. Algebra}
    \Rvolume{200}
    \Ryear{1998}
    \Rnumber{2}
    \Rpages{650-700}

\refer{Book}{StFa88}
    \Rauthor{Strade, H.; Farnsteiner, R.}
    \Rtitle{Modular {L}ie algebras and their representations}
    \Rseries{Monographs and Textbooks in Pure and Applied Mathematics}
    \Rvolume{116}
    \Rpublisher{Marcel Dekker, Inc., New York}
    \Ryear{1988}
    \Rpages{x+301}

\refer{Other}{GAP}
    \Rauthor{The GAP Group}
    \Rtitle{GAP -- Groups, Algorithms and Programming, Version 4.11.1. Available from http://www.gap-system.org}
    \Ryear{2021}

\refer{Paper}{VL06}
    \Rauthor{Vaughan-Lee, M.}
    \Rtitle{Simple Lie algebras of low dimension over GF(2)}
    \Rjournal{LMS J. Comput. Math.}
    \Rvolume{9}
    \Ryear{2006}
    \Rpages{174-192}

\refer{Paper}{ZhLin92}
    \Rauthor{Zhang, Y. Z.; Lin, L.}
    \Rtitle{Lie algebra $K(n,\mu_j,m)$ of Cartan type of characteristic $p=2$}
    \Rjournal{Chinese Ann. Math. Ser. B}
    \Rvolume{13}
    \Ryear{1992}
    \Rnumber{3}
    \Rpages{315-326}

\end{references}
